\documentclass[11pt]{article}                              %[18.06.02]

\makeatletter
\def\draft{\textheight=10.5truein \textwidth=7.5truein \parindent=8pt
           \voffset=-1truein \topmargin=0Truein
           \ifcase \@ptsize \hoffset=-1.5truein \or \hoffset=-1.35truein
                        \or \hoffset=-1.15truein \fi}
\def\quality{\textheight=240mm \textwidth=160mm \topmargin=0Truein
             \ifcase \@ptsize \hoffset=-23mm \or \hoffset=-20mm
                          \or \hoffset=-15mm \fi}
\def\USdraft{\textheight=9.85truein \textwidth=7.5truein \parindent=8pt
             \voffset=-1.0truein \topmargin=0Truein
             \ifcase \@ptsize \hoffset=-1.5truein \or \hoffset=-1.35truein
                          \or \hoffset=-1.15truein \fi}
\def\USquality{\textheight=250mm \textwidth=170mm \topmargin=0Truein
               \voffset=-1truein
               \ifcase \@ptsize \hoffset=-23mm
                       \or \hoffset=-20mm \or \hoffset=-15mm \fi}
\makeatother\draft\quality

\def\n{\noindent}  \def\const{\, {\rm Const} \,}  \def\ep{\varepsilon}
\def\IN{\hbox{I\kern-.2em\hbox{N}}}               \def\CR{$$ $$}
\def\IR{\hbox{\rm I\kern-.2em\hbox{\rm R}}}
\def\IZ{\hbox{{\rm Z}\kern-.3em{\rm Z}}}
\def\dist{\, {\rm dist}}
\def\beq#1#2{\begin{equation} \label{#1} #2 \end{equation}}
\def\bea#1{\begin{eqnarray*} #1 \end{eqnarray*}} \def\a{\!\!\!&\!\!\!\!&}
\def\function#1{\left\{\!\!\!\begin{array}{ll} #1 \end{array} \right.}
\def\proof{\smallskip \noindent {\bf Proof. \ }}       %start of proof
\def\blanksquare{\,\,\,$\sqcup\!\!\!\!\sqcap$}         %blank  square
\def\qed{\hfill\blanksquare\linebreak\smallskip\par}   %end of proof

\def\bline(#1,#2)(#3,#4)(#5){\put(#1,#2){\line(#3,#4){#5}}}  %straight line
\newcommand\mlbscale{1pt} %to change: \renewcommand\mlbscale{1.3pt}
\newif\iffigs\figstrue %\newif\iffigs\figsfalse -- to fake figures

\def\Bfig(#1,#2)#3#4{\begin{figure} \begin{center}
    \setlength{\unitlength}{\mlbscale}
       \iffigs \begin{picture}(#1,#2) #3 \end{picture}
       \else \begin{picture}(60,10)(0,0)
                   \put(0,0){\framebox(60,10){Figure}} \end{picture} \fi
    \end{center} \caption{#4} \end{figure}}
\def\bpic(#1,#2)#3{\setlength{\unitlength}{\mlbscale}
    \begin{picture}(#1,#2) #3 \end{picture}}
\def\thname{Theorem}     \def\lmname{Lemma}      \def\prname{Proposition}
\def\dfname{Definition}  \def\crname{Corollary}  \def\rmname{Remark}

\newtheorem{theorem}{\thname}[section]   %Numbering: Theorem--Other section
\newtheorem{lemma}{\lmname}[section]     %{lemma}[theorem]{Lemma}   section
\newtheorem{corollary}[lemma]{\crname}   %lemma

\newtheorem{dftn}{\dfname}[section]
\newtheorem{rmrk}[lemma]{\rmname}
             \catcode`@=11 \@addtoreset{equation}{section} \catcode`@=12

\def\c#1{#1^*} \def\cP{{\cal P}} \def\map{T}  \def\noprint#1{}
\def\Reg{{\rm Reg}} \def\Per{{\rm Per}}  \def\Free{{\rm Free}} \def\cB{{\cal B}}
\def\cA{{\cal A}}       \def\cM{{\cal M}}
\def\per#1{\langle#1\rangle}  
\def\intp#1{\left\lfloor#1\right\rfloor}
\def\toas#1{\stackrel{#1}{\longrightarrow}}
\def\blim#1#2{\if #1+ \limsup_{#2} \else {\if #1- \liminf_{#2} \else
                  \lim_{#2}\left(\begin{array}{l}\sup\\
                           \inf\end{array}\right) \fi} \fi} %\blim\pm{t\to0}

\begin{document}
%%%--------------------------------------------------------%%%%

\title{Ergodic properties of a simple deterministic \\
       traffic flow model re(al)visited}
\author{Michael Blank\thanks{
        Russian Academy of Sci., Inst. for
        Information Transm. Problems, 
        and Observatoire de la Cote d'Azur, ~ 
        e-mail: blank@iitp.ru}
        \thanks{This research has been partially 
                supported by RFFI and CRDF grants
        and a part of it has been done during my stay at ESI (May, 2002)
               .}
       }
\date{June, 2002} %\today} %
\maketitle

\n{\bf Abstract}. We study statistical properties of a family of
maps acting in the space of integer valued sequences, which model
dynamics of simple deterministic traffic flows. We obtain
asymptotic (as time goes to infinity) properties of trajectories
of those maps corresponding to arbitrary initial configurations
in terms of statistics of densities of various patterns and
describe weak attractors of these systems and the rate of
convergence to them. Previously only the so called regular
initial configurations (having a density with only finite
fluctuations of partial sums around it) in the case of a slow
particles model (with the maximal velocity 1) have been studied
rigorously. Applying ideas borrowed from substitution dynamics we
are able to reduce the analysis of the traffic flow models
corresponding to the multi-lane traffic and to the flow with fast
particles (with velocities greater than 1) to the simplest case
of the flow with the one-lane traffic and slow particles, where
the crucial technical step is the derivation of the exact life-time 
for a given cluster of particles. Applications to the optimal 
redirection of the multi-lane traffic flow are discussed as well.

\bigskip%
\n{\bf Keywords}: dynamical system, traffic flow, substitution dynamics, 
attractor, rate of convergence, large deviations.

\bigskip%
\n{\bf AMS Subject Classification}: Primary 37A60; Secondary
37B15, 37A50, 60K35.

\section{Introduction}\label{s:intro}

Let $X_M:=\{x=(\dots x_{-1} x_0 x_1\dots): \quad x_i\in\cA_M, ~
i\in\IZ\}$ be the space of be-infinite sequences (which we also
call configurations) from the alphabet $\cA_M=\{0,1,2,\dots,M\}$.
We equip this space with the metric
$$ \dist_{M}(x,y):=\sum_{i=-\infty}^\infty(M+1)^{-|i|}|x_i-y_i| $$
and consider a map $\map_{1,M}:X_M\to X_M$ from this metric space 
into itself:
\beq{def:fast-one}{
 (\map_{1,M} x)_i
     := x_i + \min\{x_{i-1}, M-x_i\} - \min\{x_i, M-x_{i+1}\} .}%
One can interpret the $i$-th coordinate of $x\in X_M$ as $x_i$
particles and $M-x_i$ holes (empty places) located at the site $i$ of 
the integer lattice $\IZ$. Then this map can be considered as a discrete
time / discrete space model for multi-lane highway traffic when
a particle (vehicle) at site $i$ of the lane $j$ can switch to any
other lane $j'$ (nonnecessary neighboring) whenever it does not disrupt
the motion of other particles, i.e. the sites $i,i+1$ of the lane $j'$
are not occupied.
From the point of view of probability theory the dynamics of this map
is a deterministic version of an asymmetric exclusion process,
i.e. the motion of a collection of random walkers constrained to the
nonintersection assumption (see e.g. \cite{Li,DLS}). Traffic flow phenomena 
have attracted considerable interest during last decade both from the 
applied and theoretical points of view. For the general account on 
these matters we refer the reader to recent reviews \cite{GG, CSS} 
(and numerous references cited there) and in this paper we shall 
concentrate only on the mathematical background of deterministic 
models of traffic flows.

We shall refer to the system $(\map_{1,M},X_M)$ as the 
{\em slow particles model}, and to take into account traffic flows 
where particles can move with the (maximal) velocity $v>1$ 
(a {\em fast particles model}) we consider a family of maps 
$\map_{v,M}:X_M\to X_M$ describing the $M$-lane traffic flow 
model with the maximal velocity $|v|$, i.e. a particle in 
this flow can move to the right (left if $v<0$) by at most $|v|$ 
positions if those positions are not occupied. 

To simplify the notation we shall drop the indices if they are 
equal to $1$, i.e. $\map_{2}$ means the case $v=2,M=1$ and 
$\map_{1,3}$ means the case $v=1,M=3$, while $\map$ stands for 
the case $v=M=1$. By a {\em dual} configuration for the 
configuration $x\in X_{M}$ we mean a configuration $\c{x}\in X_{M}$ 
such that $\c{x}_{i}=M-x_{i} ~ \forall i$. The operation of 
taking a dual can be applied also to the map by means of the relation:
$\c{\map_{v,M}}\c{x}:=\c{(\map_{v,M}x)}$ $\forall x\in X_M$.

To illustrate the usage of the dual operation consider the slow
particles model with `smart drivers', who anticipating the motion
of at most $m$ cars ahead, may move to an occupied site ahead of it
with the maximal velocity $1$. Example for the case $m=2$:
$\per{01110}\to\per{01011}$, where $\per{\cdot}$ denotes the main 
period of a (space) periodic configuration. It is straightforward 
to show that this model is described by the map 
$\tilde\map x = \c{(\map_{-m,1}\c{x})}=\c{\map_{-m,1}}x$.

By a {\em word} $A$ we shall call any (finite or infinite) sequence of
elements $a_i\in\cA_M$ and introduce the notion of the {\em density} of
a finite word $A$ in a finite word $B$ as%
\beq{def:density}{
 \rho(B,A) := \frac1{|B|}\sum_{i=1}^{|B|-|A|}
              {\rm min}_{j=1}^{|A|}
              \left\{\intp{\frac{B_{i-1+j}}{A_j}}\right\} ,}%
where $|A|$ is the length of the word $A$, $A_j\in\cA_M$ is the $j$-th
element of the word $A$, $\intp{a}$ is the integer part of the number $a$,
and we set $0/0\equiv1$ here and in the sequel. In the case $M=1$ the number
$\rho(A,B)\in[0,1]$ and is equal to the number of occurencies of
the subword $B$ in $A$ divided by the length of $A$, while in the
general case $M>1$ we have $\rho(A,B)\in[0,M]$ and the formula
(\ref{def:density}) takes into account multiplicities of those
occurencies. Example: $\rho(255,12)=\frac13(2+2)=4/3$.

The generalization of this notion for an infinite word/configuration
$x\in X_M$ leads to the notion of {\em lower/upper density}:
$$ \rho_{\pm}(x,A) := \blim{\pm}{n,m\to\infty}~\rho(x[-n,m],A), $$
where (and in the sequel) $\limsup$ corresponds to the index $+$
and $\liminf$ to the index $-$, and $x[n,m]$ a subword of the word $x$
which starts from the position $n$ and goes till the position $m$ in the
original word.
%$$ \rho_-(A,B) := \liminf_{n,m\to\infty}\rho(A[-n,m],B), \qquad
%   \rho_+(A,B) := \limsup_{n,m\to\infty}\rho(A[-n,m],B) .$$
The asymmetry with respect to $n$ and $m$ is necessary to take into account
the possibility to have left and right `tails' with different statistics:
for $x:=\dots 00001111\dots$ we have $\rho_-(x,1)=0$ and $\rho_+(x,1)=1$,
while $\rho(x[-n,n],B)\toas{n\to\infty}1/2$.
Observe also that for a (space) periodic configuration
$\per{A}:=\dots AAA\dots$ we have
$\rho_-(\per{A},B)=\rho_+(\per{A},B)\equiv\rho(\per{A},B)$ for any
pair of finite words $A,B$.

For a system of particles on a lattice one can define its average velocity
as follows. First, for each particle in a configuration $x\in X$ we define
its `local' velocity as a distance by which it will move on the next
step of the dynamics:
$$ V(x,i) := \min\{v, ~ x_i, ~ \min\{j-i: ~ j>i, x_j=1\}\} ,$$
and, since for $M>1$ a site $i$ in the configuration $x\in X_M$
may contain several particles (i.e. $x_i>1$), we sum up their
velocities to get $V(x,i)$. For example, in the case of the map
$\map_{1,M}$ we have $V(x,i) := \min\{x_i, M-x_{i+1}\}$. Note
that the `local' velocity is well defined for each site $i$ of a
configuration $x\in X_M$ (independently on the presense of a
particle there), indeed, if $x_i=0$ (i.e. there is no particle at
this site) we have $V(x,i)=0$. Now we define the lower/upper {\em
average velocity} as
$$ V_{\pm}(x):=\blim{\pm}{n,m\to\infty}~
               \frac1{\rho(x[-n,m],1)\cdot(n+m+1)}\sum_{i=-n}^m V(x,i) .$$
%$$ V_-(x):=\liminf_{n,m\to\infty}
%           \frac1{\rho(x[-n,m],1)\cdot(n+m+1)}\sum_{i=-n}^m V(x,i), \CR
%   V_+(x):=\limsup_{n,m\to\infty}
%           \frac1{\rho(x[-n,m],1)\cdot(n+m+1)}\sum_{i=-n}^m V(x,i) .$$
Often it is more suitable to work with another statistics, called
{\em flux}, equal to the number of particles crossing a given
position on the lattice per unit time, i.e.
$\Phi(x[-n,m]):=\frac1{n+m+1}\sum_{i=-n}^m V(x,i)$.
Thus we define the upper/lower {\em average flux} as
$$ \Phi_{\pm}(x):=\blim{\pm}{n,m\to\infty}~\frac1{n+m+1}\sum_{i=-n}^m V(x,i) .$$
%$$ \Phi_-(x):=\liminf_{n,m\to\infty}\frac1{n+m+1} \sum_{i=-n}^m V(x,i), \quad
%   \Phi_+(x):=\limsup_{n,m\to\infty}\frac1{n+m+1} \sum_{i=-n}^m V(x,i) .$$%
We shall use also the notation $\Phi_\pm^{(v)}$ to indicate the
maximum velocity if needed, and $0_i:=\underbrace{00\dots0}_{i}$. 
The connection of the flux to the densities is given by the following 
simple result. 

\begin{lemma}
$\Phi_\pm^{(v)}(x)=\sum_{i=1}^v\rho_\pm(x,10_i)$ for $x\in X$,
in particular $\Phi_\pm^{(1)}(x):=\rho_\pm(x,10)$.
\end{lemma}

\proof By definition we have%
\bea{ \Phi_\pm^{(v)}(x)
 \a= v\cdot\rho_\pm(x,10_v) + (v-1)\cdot(\rho_\pm(x,10_{v-1})-\rho_\pm(x,10_v))
                          + \dots + 1\cdot\rho_\pm(x,10) \\
 \a= (v-(v-1))\cdot\rho_\pm(x,10_v) + ((v-1)-(v-2))\rho_\pm(x,10_{v-1})
                          + \dots + 1\cdot\rho_\pm(x,10) \\
 \a= \sum_{i=1}^v\rho_\pm(x,10_i) .}%
\qed

%%%%%%%%%%%%%%%%%%%%%%%%%%%%%%%%%%%%%%%%%%%%%%%%%%%%%%%%%%%%%%%%%%
%% Fundamental diagram
\Bfig(320,100)
      {\footnotesize{
       %% picture (a)
       \bline(0,0)(1,0)(320)   \bline(0,0)(0,1)(100)
       \bline(0,100)(1,0)(320) \bline(160,0)(0,1)(100) \bline(320,0)(0,1)(100)
       \put(20,20){\vector(1,0){120}} \put(20,20){\vector(0,1){70}}
       \bline(20,75)(1,0)(15)  \bezier{100}(35,75)(70,23)(130,20)
       \put(145,18){$\rho$}    \put(28,10){$\frac{M}{v+1}$}  \put(123,12){M}
       \put(12,85){$V$}        \put(5,72){$Mv$}  \put(75,45){$\frac{M}\rho-M$}
       \bezier{30}(35,75)(35,50)(35,20)
       % picture (b)
       \put(160,0){\bpic(150,100){
       \put(20,20){\vector(1,0){120}} \put(20,20){\vector(0,1){70}}
       \bline(20,20)(1,4)(13.7)  \bline(125,20)(-5,3)(91)
       \put(145,18){$\rho$}   \put(27,10){$\frac{M}{v+1}$}  \put(123,12){M}
       \put(12,85){$\Phi$}    \put(2,71){$\frac{Mv}{v+1}$} \put(75,53){$M-\rho$}
       \put(20,57){$v\rho$}
       \bezier{30}(33.5,74)(33.5,50)(33.5,20)
       \bezier{10}(33.5,74.5)(27,74.5)(20,74.5)
       }}
      }}
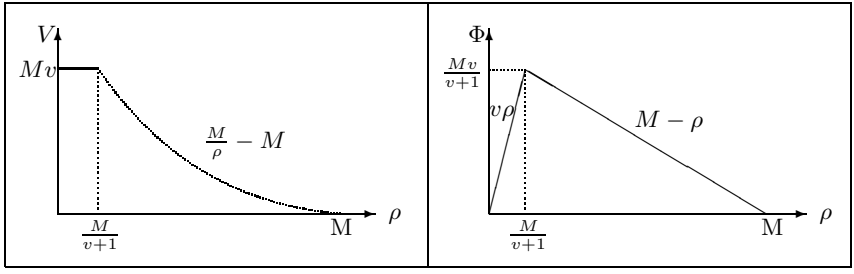
{Fundamental diagram for $\map_{v,M}$: dependence of the average
velocity $V$ or the flux $\Phi$ on the density
$\rho=\rho(\cdot,1)$. \label{av-speed-den}}
%%%%%%%%%%%%%%%%%%%%%%%%%%%%%%%%%%%%%%%%%%%%%%%%%%%%%%%%%%%%%%%%%

The main results of the paper are the following statements.

\begin{theorem}\label{t:density-pres} (Invariance of densities)
$\rho_\pm(\map_{v,M}^tx,A)=\rho_\pm(x,A)$ for all $x\in X_M$ and 
$t\in\IZ_+$ if and only if $A\in\{0,1\}$.
\end{theorem}

Denote by $\Free_v:=\{x\in X_{M}: ~ V(x,i)=v\cdot x_i ~ \forall
i\in\IZ\}$ the subset of configurations where all particles have
the maximal available velocity and thus move independently.
Clearly, $\map_{v.M}(\Free_v)=\Free_v$ and
$\map_{v.M}(\c{\Free_v})=\c{\Free_v}$. It is of interest that
$\Free\cap\c{\Free}\ne\emptyset$.

\begin{theorem}\label{t:conf-conv} (Convergence)
The set \/ $\Free_v\cup\c{\Free}_v$ is the only locally maximal weak
attractor of the dynamical system $(\map_{v,M},X_M)$, and for $x\in X_M$
we have
$\map_{v,M}^tx \toas{t\to\infty}
      \function{\Free_v     &\mbox{if } \rho_+(x,1)\le\frac{M}{v+1} \\
                \c{\Free}_v &\mbox{if } \rho_-(x,1)\ge\frac{M}{v+1} }$.
\end{theorem}

\begin{theorem}\label{t:flux-conv} (Limit flux)
$\Phi_\pm^{(v)}(\map_{v,M}^tx)\toas{t\to\infty}F_{v,M}(\rho_\pm(x,1))$,
where
$F_{v,M}(\xi):=\function{v\xi  &\mbox{if } \xi\le \frac{M}{v+1} \\
                         M-\xi &\mbox{otherwise}}$.
\end{theorem}

Denote by $\mu_p$ a product (Bernoulli) measure with the density 
$pM$ on the space of sequences $X_M$.

\begin{theorem}\label{t:typ-conv} (Typical dynamics)
For $\mu_p$-a.a. $x\in X_M$ we have $\rho(x,1)=pM$ and 
$\dist_M(\map_{v,M}^tx, \Free\cup\c\Free) \le M^{-t/\gamma+1}$
and
$\blim+{n\to\infty}\frac1{2n}\sum_{i=-n}^nV(\map_{v,M}^nx,1)
  = F_{v,M}(\rho(x,1))$ for any $\gamma\in(0,1)$.
\end{theorem}

In this result one can use instead of $\mu_p$ any probabilistic 
translation invariant measure with fast enough decay of 
correlations (see Lemma~\ref{l:gen-mes-density}).

Proofs of Theorems~\ref{t:density-pres}--\ref{t:typ-conv} are based 
on the reduction of the general case $v,M\ge1$ to the simplest one 
$v=M=1$. For $v=1,~M>1$ this reduction boils down to the proof that 
a multi lane traffic flow can be represented by a direct product 
of one-lane flows (see Theorem~\ref{t:sawtooth} describing the 
`sawtooth redirection' construction). In the case $v>1, M=1$ we 
make use of a specially constructed substitution dynamics 
(see Lemma~\ref{l:fast-equiv}) to prove the reduction, while in the 
general case $v,M>1$ we combine these two arguments. The main technical 
step of the analysis in the case $v=M=1$ is the derivation of the 
exact life-time for a given cluster of particles, i.e. the number 
of iterations after which it will disappear, described in 
Lemma~\ref{l:life-time}.
Note that earlier only very weak (and unnaturally large) estimates
of the life-time type were known (see, e.g. \cite{Bl-var,Bl-jam,Fu}).

We provide also the analysis of the rate of convergence to the limit
of various statistics for space periodic, regular, and `typical' 
initial configurations based on large deviations estimates 
(see Section~\ref{s:rate}), and study the dynamics of a passive 
tracer in the flow of fast particles (Section~\ref{s:tracer}).

It is clear that the main problem in the study of traffic flows
is the analysis of `traffic jams' (without them the dynamics is
trivial): we shall say that a segment $x[n,m]$ with $m>n$
corresponds to the {\em jammed cluster} if
$\max\{V(x,i)/x_{i}: ~ n\le i< m\}<v$ and
$V(x,n-v)/x_{n-v}=V(x,m)/x_{m}=v$, i.e. all particles inside of
this segment do not have the maximum available velocity. Note
that in the case $v=M=1$ the jammed cluster is the same as the
cluster of particles.

\section{Dynamics of slow particles $(\map,X)$}\label{s:slow}

This model has been introduced originally in \cite{NS} for the case
of a traffic flow on a finite lattice (say of size $L$) with periodic
boundary conditions and studied numerically in a large number of
publications. It is straightforward to show that this case
corresponds to the restriction of the map $\map$ to (space) $L$-periodic
configurations. The first `quasi'-analytic result for the $L$-periodic
case has been obtained in \cite{Fu} for `typical' initial
configurations of length $L$. However the first complete proof
appeared only in \cite{Bl-var}, where regular initial configurations
on the infinite lattice were considered as well. In this section
we shall study the problem for all initial configurations, using
a rather different and more simple approach than the one in \cite{Bl-var}.

Let us start from the analysis of lower and upper densities. Note
that if the lower density coincides with the upper one, i.e. the 
limit value exists, we cal this common value the density 
$\rho(\cdot,\cdot)$. Example when they do not coincide: 
$1=\rho_{+}(\dots 111000\dots,1) \ne \rho_{-}(\dots111000\dots,1)=0$.

\begin{lemma}\label{l:dens-ineq} $\rho_\pm(x,1)=1-\rho_\mp(x,0)$,
and thus $\rho_\pm(\c{x},1)=1-\rho_\mp(x,1)$.\end{lemma}

\proof By the definition of the lower density we have
$$ \rho_-(x,1) = \blim-{n,m\to\infty}\rho(x[-n,m],1)
 = 1 - \blim+{n,m\to\infty}\rho(x[-n,m],0)
 = 1 - \rho_+(x,0) ,$$
since $\rho(A,1)\cdot|A| + \rho(A,0)\cdot|A| = |A|$ for any finite
binary word $A$. The derivation for the upper density follows the
same argument, while the second statement follows from the identity:
$\rho(\c{x}[-n,m],1)=\rho(x[-n,m],0)=1-\rho(x[-n,m],1)$. \qed

\begin{lemma}\label{l:enclos-ineq}
$\rho_\pm(x,A)\ge\rho_\pm(x,B)\cdot\rho(B,A)$ for any configuration 
$x\in X$ and any pair of finite words $A,B$. \end{lemma}
\proof If $A\not\subseteq B$ the inequality becomes trivial, since
$\rho(B,A)=0$, while $\rho_\pm(x,A)\ge0$.
Assume now that $A\subseteq B$. Then $\forall n,m\in\IZ$ we
have $\rho(x[-n,m],A)\ge\rho(x[-n,m],B)\cdot\rho(B,A)$ because the right
hand side takes into account only those enclosures of $B$ to $x$ when
the word $B$ belongs to a segment $x[i,j]=A$, while there might be
other enclosures as well. \qed

\n{\bf Proof} of Theorem~\ref{t:density-pres} in the case $v=M=1$. 
Let us prove first that $\rho_{\pm}(x,1)=\rho_{\pm}(\map{x},1)$
for all $x\in X$. For any $n,m\in\IZ_{+}$ we have
$|\sum_{i=-n}^{m}(x_{i}-(\map x)_{i})|\le 2$, since during one
iteration of the map at most one particle can enter the interval
of sites from $-n$ to $m$ (from behind) and at most one particle
can leave this interval.

By the definition of the lower density there is a sequence of pairs
$(n_{j},m_{j})\toas{j\to\infty}(\infty,\infty)$ such that
$$ \frac1{n_{j}+m_{j}+1}\sum_{i=-n_{j}}^{m_{j}}x_{i}
   \toas{j\to\infty} \rho_{-}(x,1) .$$
On the other hand, since
$|\sum_{i=-n_{j}}^{m_{j}}x_{i} - \sum_{i=-n_{j}}^{m_{j}}(\map x)_{i}| \le2$,
we deduce that $\rho_{-}(x,1)$ is a limit point for partial sums for
the sequence $\map x$. Therefore we need to show only that this is
indeed the lower limit. Assume, on the contrary, that there is
another limit point, call it $\xi$, for the partial sums for $\map x$ 
such that $\xi<\rho_{-}(x,1)$. Doing the same operations with the
partial sums for $\map x$ converging to $\xi$ we can show that this 
value is also a limit point for the partial sums for the sequence
$x$, and, hence, $\xi$ cannot be smaller than $\rho_{-}(x,1)$.

The proof for the upper density follows from the similar argument.

By Lemma~\ref{l:dens-ineq} we have $\rho_\pm(x,1)=1-\rho_\mp(x,0)$, 
which proves the preservation of the density of zeros as well.

To prove that all other statistics are not preserved under dynamics
we need to study it in more detail. Therefore we postpone
the continuation of the proof till the end of this Section.

\begin{lemma}\label{l:dual-slow} $\c{\map}=\map_{-1}$.
\end{lemma}

\proof The action of the map $\map$ on binary configuration is
equivalent to the exchange of any pair $10$ to $01$. Since the dual
map describes the dynamics of holes it corresponds in this case to
the exchange of pairs $01$ to $10$, which proves the statement. \qed

By a {\em cluster} (of particles) in a binary configuration
$x\in X$ we mean a collection of consecutive positions
$x_{m},x_{m+1},\dots,x_{n}$ such that $n-m>1$, $x_{i}=0$ $\forall
i\in(m,\dots,n)$ and $x_{m-1}=x_{n+1}=0$. After each iteration of
the map $\map$ the last particle in the cluster moves away (i.e.,
$(\map x)_{n}=0$) and either appears a new element in the cluster
from the left (i.e., $(\map x)_{m-1}=1$ and $(\map x)_{m-2}=0$),
or the first particle preserves its position $m$. Therefore the 
number of particles in a given cluster cannot increase, and the time
up to the moment when the cluster length shrinks to 1 (i.e. it 
disappears) we shall call the {\em life-time} of the cluster or the 
number of iterations which are needed for a given cluster to disappear.

Define an integer-valued function %
\beq{l-count}{ I(x,i):=\max\{k<i: ~ \rho(x[k,i],1)=\rho(x[k,i],0) \} .}

\begin{lemma}\label{l:life-time}
Let $\sum_{i=m}^{n}x_i=n-m+1$ (i.e. the positions from $m$ to $n$
correspond to a cluster of particles and let $I(x,n)>-\infty$,
then after exactly $\frac12(n-I(x,n)-1)$ iterations (which is equal to
the number of ones minus one in the word $x[n-I(x,n),n]$) this cluster
will disappear. If $\rho_{+}(x)\le1/2$ then $\forall i\in\IZ$ we have
$I(x,i)>-\infty$. \end{lemma}

\proof Let
$\Omega^{2n} := \{A\in\{0,1\}^{2n}: ~ A_{2n}=1, ~ I(A,2n)=1 \}$,
where $I(x,i)$ is defined by the relation~(\ref{l-count}) and thus
$\rho(A,1)=|A|/2$.
Observe that if $A\in\Omega^{2n}$ than for any $0<m<n$ and any word
$B\in\{0,1\}^{2m}$ such that $B_i=A_i$ $\forall 0<i\le 2m$ we have
$B\not\in\Omega^{2m}$. Therefore we shall call the words from
$\Omega^{2n}$ {\em minimal} words (or minimal intervals) corresponding
to clusters of particles in their ends.

Consider a map $\Gamma:\Omega^{2n}\to\IZ^{2n-2}$ defined by the relation:
$$ (\Gamma A)_{i+1} = A_i + \min\{A_{i-1}, 1-A_i\} - \min\{A_i, 1-A_{i+1}\} .$$
Observe that this is a shift to the right of the action of our map $\map$.
We shall prove that for each $n$ we have
$\Gamma:\Omega^{2n}\to\Omega^{2n-2}$.
%Examples of the action of $\Gamma$:%
%\begin{verbatim}
%0001011011      00100111
% 00011011        001011
%  001011          0011
%   0011            01
%    01
%\end{verbatim}.

\bigskip

\begin{center}
\begin{tabular}{rlllllllllllllllllll}
$A=$        &0&0&0&1&0&1&1&0&1&1& \quad &0&0&1&0&0&1&1&1 \\
$\Gamma^1A=$& &0&0&0&1&1&0&1&1& &       & &0&0&1&0&1&1&  \\
$\Gamma^2A=$& & &0&0&1&0&1&1& & &       & & &0&0&1&1&    \\
$\Gamma^3A=$& & & &0&0&1&1& & & &       & & & &0&1&      \\
$\Gamma^4A=$& & & & &0&1& & & & &                        \\
\end{tabular}\par%
Examples of the action of $\Gamma^t$ on $\Omega^{10}$ (left) and
$\Omega^{8}$ (right).
\end{center}\bigskip

Let $A\in\Omega^{2n}$ and let $\zeta_A$ be the position of the last $0$ in $A$.
For each word $A\in\Omega^{2n}$ define a new word $A'\in\IZ^{2n-2}$ as follows:
$$ A'_i := \function{A_i &\mbox{if } i\le 2n-2 \mbox{  and } i\ne\zeta_A \\
                       1 &\mbox{if } i=\zeta_A .} $$
Then $|A'|=|A|-2=2n-2$, $\rho(A',1)=1/2$ and thus $A'\in\Omega^{2n-2}$
since otherwise $A$ would be not minimal as well. Now using the following
simple identity:
$$ (\Gamma A)_i=\function{(\Gamma A')_i &\mbox{if } \zeta_A-i\not\in\{1,2\} \\
                                 0      &\mbox{if } i=\zeta_A-2 \\
                            A_{\zeta_A} &\mbox{if } i=\zeta_A-1 \\
                                 1      &\mbox{if } i=|A|-2 ,} $$
we get $\Gamma A\in\Omega^{2n-2}$. Example: $A=001011$, $A'=0011$.

It remains to show that if we have a cluster of particles located in
the end of a minimal configuration $A\in\Omega^{2n}$ then this cluster 
(i.e. particles at sites from $\zeta_A+1$ to $2n$) will vanish after
$n-1$ iterations. Observe that after one iteration of the map $\map$
the cluster either preserves its length, or the length decreases by
one (when two positions immediately preceeding the cluster are occupied by
two zeros). The map $\Gamma$ defined above controls this process since for
each $A\in\Omega^{2n}$ and each $0\le t\le n-1$ the last positions starting
from $(\zeta_{\Gamma^tA}+1)$ correspond to the cluster under study. \qed

\begin{corollary} The dynamics of a cluster of particles depends only 
on the distribution of particles lying below the cluster. Moreover, for 
a given cluster of particles only particles belonging to its minimal 
word can join the cluster. Thus, if a particle does not belong to any 
minimal word, then for each $t\in\IZ_+\cup\{0\}$ its local velocity in 
the configuration $\map^tx$ is equal to 1, i.e. it moves freely.
\end{corollary}

\begin{lemma}\label{l:min-int} Let $A:=x[n,m]$ and $A':=x[n',m']$ be
two minimal words in the configuration $x\in X$. Then the inequality
$m<m'$ yields either $m<n'$ (i.e. $A\cap A'=\emptyset$), or
$n'<n$ (i.e. $A\subset A'$).
\end{lemma}

\proof Assume on the contrary that $n\le n'\le m$. Then by the definition
of a minimal word we have%
\bea{\frac12(m'-n'+1)
 \a= \rho(x[n',m'],1)\cdot(m'-n'+1) \\
 \a= \rho(x[n',m],1)\cdot(m-n'+1) + \rho(x[m,m'],1)\cdot(m'-m+1) \\
 \a> \frac12(m-n'+1+m'-m+1) = \frac12(m'-n'+2) .}%
We came to a contradiction. \qed

\begin{corollary} For a given particle $\xi$ in a configuration $x\in X$
let the length of the largest minimal word to which $\xi$ belongs is $2n$.
Then for any $t\ge n-1$ the local velocity of the particle $\xi$ in the
configuration $\map^tx$ is equal to 1.
\end{corollary}

Observe that in the case $v=M=1$ the set $\Free=\{x\in X_1: ~
x_ix_{i+1}=0 ~ \forall i\}$ and is the union of `free' particles
(i.e. particles having velocity $1$), while its dual
$\c{\Free}=\{x\in X_1: ~ (1-x_i)(1-x_{i+1})=0 ~ \forall i\}$
corresponds to `free' holes (i.e. to holes having velocity $-1$).

\begin{lemma}\label{l:conv-low} Let $\rho_{+}(x,1)\le1/2$, then each
cluster of particles in the configuration $x$ will disappear
after a finite number of iterations and
$\dist(\map^tx,\Free)\toas{t\to\infty}0$. If there is a cluster
of particles having an infinite minimal word (i.e. which does not
vanish in finite time), then $\rho_{+}(x,1)>1/2$.
If $\rho_{+}(x,1)>1/2$ then there are clusters of particles with
arbitrary large (but may be finite) minimal words.
\end{lemma}

\proof Let $x\in X$ satisfies the assumption that $\rho_{+}(x,1)\le1/2$
and let the segment $x[n,m]$ be a cluster of particles. Then there
exist a pair of integers $n',m'$ such that $\-infty<n'<n<m\le m'<\infty$
and $x[n',m']$ is the largest minimal word covering the cluster of
particles $x[n,m]$ (otherwise this would contradict to the definition
of the upper density). Hence by Lemma~\ref{l:life-time} after at most
$(m'-n')/2$ iterations this cluster will disappear and all particles
will become free. Since this argument can be applied to any cluster of
particles, this yields the first statement.

Assume now that the minimal word of a cluster of particles $x[n,m]$
is not bounded. Then for any $n'<n$ we have $x[n',m]>1/2$ and thus
for any $k\in\IZ_+$ we have%
\bea{ \rho(x[n-k^2,m+k],1)
 \a\ge \frac{m+k-n+k^2+1}{m-n+k^2+1}\rho(x[n-k^2,m],1) \\
   \a> \frac{m-n+k^2+1}{m+k-n+k^2+1}\frac12
     = \left(1 - \frac{k}{m+k-n+k^2+1}\right)\frac12
       \toas{k\to\infty}\frac12 .}%
Therefore $\rho_{+}(x,1)>1/2$, which proves the second statement.

The last statement is an immediate consequence of the definition
of the minimal word. \qed

\begin{lemma}\label{l:conv-high} Let $\rho_{-}(x,1)>1/2$ then 
$\dist({\c\map}^t\c{x},\Free)\toas{t\to\infty}0$. \end{lemma}

\proof By Lemma~\ref{l:dens-ineq} we have
$\rho_{+}(\c{x},1)=1-\rho_{-}(x,1)<1/2$. On the other hand,
Lemma~\ref{l:dual-slow} shows that asymptotic properties of the
maps $\map$ and $\c\map$ coincide, thus we can apply the statements
of Lemma~\ref{l:conv-low} for the case of $(\c\map)^t\c{x}$ to prove
the desired result. \qed

Note now that there are configurations not satisfying the
assumptions of Lemmata~\ref{l:conv-high},\ref{l:conv-low} which
still converge to $\Free\cup\c\Free$ under the action of the map
$\map$. Indeed, let $y=\dots111000\dots$ and let the index $0$
correspond to the first $0$ in $y$. Observe that
$\rho_{-}(y,1)=0<\rho_{+}(y,1)=1$, however%
$\dist(\map^ty,\Free) = 2^{-(t+1)}(1+2^{-2}+2^{-4}+\dots)
                      = \frac23\cdot 2^{-t} \toas{t\to\infty}0$,
since for large $t$ the `central' part of $\map^ty$ will be
occupied only by free particles.
On the other hand, for $\c{y}=\dots000111\dots$ we have
$\dist(\map^t\c{y},\Free)=2^{-1}+2^{-3}+\dots+2^{-2n+1}+\dots=1/3$
for each $t\in\IZ^+$, while%
$\dist({\c\map}^t\c{y},\Free) = \frac23\cdot 2^{-t} \toas{t\to\infty}0$.

\begin{lemma}\label{l:basin} For any $x\in X$ we have
$\dist(\map^tx,\Free)\toas{t\to\infty}0$ if and only if
$\rho_+(x[1,\infty],1)\le1/2$. If additionally there exist a pair 
$n,m\in\IZ_+$ such that $\rho(x[i,i+m-1],1)\le1/2$ for each $i\ge n$ 
we have $\dist(\map^tx,\Free) \le\const 2^{-t}$.
\end{lemma}

\proof Observe that $\rho_+(x[1,\infty],1)\le1/2$ implies that there
exists $N\in\IZ$ such that the life-time for each cluster of particles
lying to the right from $N$ is finite. On the other hand, the distance
from the most right cluster of particles in $\map^tx$ located to the
left of $N$ to the position $N$ grows with $t$ linearly. This proves
the first statement and shows that the rate of convergence might be
smaller than $2^{-t}$ only if the life-time of clusters of particles
lying to the right from a sufficiently large position $N$ is not bounded.
The additional assumption guaranteers that this cannot happen, which
yields the second statement. \qed

Let $\cM(X)$ be the set of probabilistic translation invariant 
measures on $X$ and let $\mu[\phi(x)]:=\int\phi(x)~d\mu(x)$ for 
$\mu\in\cM(X)$, in particular, $\mu[x_0]:=\mu(x\in X: ~ x_0=1)$. 
Consider a subset of $\cM(X)$ corresponding to measures in the space 
of sequences with weak dependence between coordinates (exponentially 
fast decay of correlations):
$$ \cM_p(X) := \{\mu\in\cM(X): ~ \mu[x_0]=p, ~ 
   \left| \mu[x_0\cdot x_k] - \mu^2[x_0] \right|
   \le Ce^{-\alpha|k|} \} .$$
for some $C,\alpha$ and $\forall k\in\IZ$. Note that, e.g., a 
product (Bernoulli) measure $\mu_p\in\cM_p(X)$.

\begin{lemma}\label{l:gen-mes-density} 
For any $\mu\in\cM_p(X)$ we have $\rho(x,1)=p$ for $\mu$-a.a. 
$x\in X$, and thus $\mu(x\in X: ~ \rho_{-}(x)<1/2<\rho_{+}(x))=0$.
\end{lemma}

\proof Let $S_{n,m}(x):=\sum_{i=-n}^mx_i$. Then 
$\mu[S_{n,m}(x)]=p$ and by Chebyshev inequality $\forall\ep>0$ 
we have
$$ \mu(x\in X: ~ |S_{n,m}(x) - p| \ge \ep) 
   \le \frac1{\ep^2} \cdot \mu[(S_{n,m}(x) - p)^2] .$$
A straightforward calculation shows that 
$\mu[(S_{n,m}(x) - p)^2]\le \frac{C_1}{n+m+1}$ and thus
$$ \mu\left(x\in X: ~ 
   \left|\frac1{n+m+1}\sum_{i=-n}^mx_i - p \right| \ge\ep\right)
   \le \frac{C_1}{(n+m+1)\ep^2} \toas{n,m\to\infty}0 .$$

Therefore $\rho_\pm(x[-n,m],1)\toas{n,m\to\infty}p$ in probability,
which yields the existence of the density $\rho(x,1)=p$ for $\mu$-a.a.
$x\in X$ and thus the statement under question. \qed

Denote by $\Per_n(\map):=\{x\in X: ~ \map^nx=x\}$ the set of $n$-periodic
(in time) trajectories of the map $\map$ and by
$\cB(Y):=\cup_{n\ge0}\map^{-n}Y$ the {\em basin of attraction} of a
subset $Y\subset X$.

\begin{lemma}\label{l:per-str} $\map:X\to X$ is a Lipschitz continuous
map in the topology induced by the metrics $\dist(\cdot,\cdot)$.
For each $n\in\IZ_+$ there exists an $n$-periodic trajectory, and all
periodic trajectories are unstable.
${\rm Clos}(\cB(\Free \cup \c\Free))=X$,
$(\Free\cup\c\Free)\cap\Per_1(\map)=\emptyset$, and
$\mu_p(\cB(\Per_1(\map)))=0$ while ${\rm
Clos}(\cB(\Per_1(\map)))=X$.
\end{lemma}

\proof Let us start with the Lipschitz continuity. Consider two
configurations $x\ne y\in X$ and assume that $-n<0$ is the largest negative
index and $m\ge0$ is the smallest nonnegative index of sites, where they
differ, i.e. for all $-n<i<m$ we have $x_i=y_i$. Then we have
$$ 2^{-n}+2^{-m} \le \dist(x,y) \le 2(2^{-n}+2^{-m}) .$$
On the other hand, a straightforward calculation shows that the
closest to the origin from the left side differing coordinates of the
configurations $\map x$ and $\map y$ belong to the set $\{-(n+1),-n,-(n-1)\}$,
while the closest from the right side belong to $\{m-1,m,m+1\}$. Thus
$$ 2^{-(n+1)}+2^{-(m+1)} \le \dist(\map x, \map y)
                         \le 2(2^{-(n-1)}+2^{-(m-1)}) .$$
Therefore
$$ \frac14 = \frac{2^{-(n+1)}+2^{-(m+1)}}{2(2^{-n}+2^{-m})}
 \le \frac{\dist(\map x, \map y)}{\dist(x,y)}
 \le \frac{2(2^{-(n-1)}+2^{-(m-1)})}{2^{-n}+2^{-m}} = 4 .$$

For a given $n\in\IZ_+$ consider a space-periodic configuration $x\in X$
with the (space) period $n$, e.g. $x_i=x_{i+n}$ $\forall i$. Then it is
immediate to show that for any $t\in\IZ$ the configuration $\map x$ is again
space periodic with the same period $n$ and converges either to $\Free$, or
to $\c{\Free}$, depending on its density. This gives a construction of the
$n$-periodic (in time) trajectories.

The structure on the set of fixed points $\Per_1(\map)$ is a bit more
involved:
$$ \Per_1(\map):=\left\{x^{(n)}\in X: ~~
       x_i^{(n)}=\function{0 &\mbox{if }  i<n \\
                           1 &\mbox{otherwise}} \right\}.$$
Indeed, assume that $\map x=x$, then either $x$ does not have zero 
coordinates, or all coordinates starting from, say, $n$-th, should be 
equal to one. Now for $x^{(n)}\in\Per_1(\map)$ we define 
$y^{(n,m)}\in X$ such that
$$ y^{(n,m)}_i=\function{0 &\mbox{if }  i<n, \mbox{ or } i>m \\
                         1 &\mbox{otherwise}} $$ 
for some $m>n$. Then $\dist(x^{(n)},y^{(n,m)})=2^{-m}\toas{m\to\infty}0$, 
while
$\dist(\map^t x^{(n)}, \map^t y^{(n,m)})\toas{t\to\infty}2^{-(n-1)}\ne0$.
Thus for each $\ep>0$ there is configuration $x'=x'(\ep)$ such that
$\dist(x^{(n)},x')\le\ep$ and $\map^t x' \not\to x^{(n)}$ as $t\to\infty$,
which yields instability.

Observe now that the set
$Y:=X\setminus(\cB(\Free \cup \c\Free))
   = \{x\in X: ~ \rho_-(x)<1/2<\rho_+(x)\}$ has $\mu_p$-measure zero,
since for each $\mu_p$-typical trajectory the lower and upper densities
coincide. %

Consider now an arbitrary configuration $x\in X$ and a sequence
of configurations $\{y^{(n)}\}_n$ defined as
$y^{(n)}_i=\function{x_i &\mbox{if }  i<n, \\
                     1   &\mbox{otherwise}}$.
Then $\dist(x,y^{(n)})\le 2^{-n+1}$, on the other hand,
$y^{(n)}\toas{n\to\infty}\Per_1(\map)$, which proves the last statement. \qed

\n{\bf Remark}. In the case of a finite cluster of particles its
last particle immediately leaves under the dynamics. This is not the case
for clusters not bounded from the right, which explains the existence
of fixed points.

\smallskip

We shall say that a closed $\map$-invariant set $Y$ is a 
{\em weak attractor} if $\mu_{{\rm ref}}(\cB(Y))>0$.
A weak attractor $Y$ is called a {\em Milnor} attractor if
$\mu_{{\rm ref}}(\cB(Y)\setminus\cB(Y'))>0$ for any proper
compact invariant subset $Y'\subseteq Y$ (see, e.g. \cite{Mi}).

\begin{lemma} The set $\Free\cup\c{\Free}$ is an week attractor with
respect $\mu_{{\rm ref}}=\mu_p$, but not a Milnor one, moreover it
is not a topological attractor.
\end{lemma}

\proof The sets $\Free$ and $\c\Free$ are closed, since they
contain all their limit points. Let $Z^{(p)}:=\{x\in \Free: ~
\rho(x)=p\}$ -- this is an invariant set and
$\mu_p(\cB(Z^{(p)}))=1$. Denote now by $Z'^{(p)}$ a single
configuration from $Z^{(p)}$ together with all its left and right
shifts. Clearly $\mu_p(Z^{(p)}\setminus(Z^{(p)}\setminus
Z'^{(p)}))=\mu_p(Z'^{(p)})=0$. Observe that the points from the
complement to the basins of attraction of $\Free$ and $\c{\Free}$
are everywhere dense, which proves the absense of included open
sets. The last statement follows from the fact that the basin of
attraction does not contain any open set. \qed

\n{\bf Proof} of Theorem~\ref{t:density-pres} (continuation).
Let us prove now that for any word $A$ with $|A|>1$ the density
$\rho_\pm(x,A)$ is not preserved under dynamics. There might be 3
possibilities: $\rho(A,1)<1/2$, $\rho(A,1)>1/2$ and $\rho(A,1)=1/2$.
We start from the first case. Clearly $\rho(A,1)<1/2$ yields $\rho(A,00)>0$.
Consider a configuration $x:=\per{A\underbrace{11\dots1}_{2|A|}}$,
where $x=\per{B}\equiv\dots BBB\dots$ stays for a space-periodic
configuration. By the construction
$\rho(x,1)=(\rho(A,1)\cdot|A|+2|A|)/(3|A|)\ge2/3>1/2$.
Therefore $\map^t x\toas{t\to\infty}\c{\Free}$ and hence
$\rho(\map^t x,00)\toas{t\to\infty}0$. Assume now that the density is
preserved, i.e. $\rho(\map^t x,A)=\rho(x,A)$ $\forall t$. Then by
Lemma~\ref{l:enclos-ineq} we have
$$ \rho(\map^t x,00) \ge \rho(\map^t x,A)\cdot\rho(A,00)
                       = \rho(x,A)\cdot\rho(A,00) > 0 ,$$
while the left hand side vanishes when $t\to\infty$. We came to the
contradiction.

If $\rho(A,1)>1/2$ we shall follow a similar argument, considering
another space-periodic configuration 
$x:=\per{A\underbrace{00\dots0}_{2|A|}}$.

In a more delicate case $\rho(A,1)=1/2$ we do the folowing. If
additionally $\rho(A,11)>0$ we follow the same argument as in the
case $\rho(A,1)<1/2$ to show that $\rho(x,11)>0$, while
$\rho(\map^tx,11)\toas{t\to\infty}0$. If $\rho(A,00)>0$ we follow
the case $\rho(A,1)>1/2$ to show that $\rho(x,00)>0$, while
$\rho(\map^x,00)\toas{t\to\infty}0$. It remains to consider the
case when $\rho(x,11)=\rho(x,00)=0$, i.e. $A=1010\dots10$ or
$A=0101\dots01$. In the first of these case we choose
$x:=\per{1A0}$. Then
$$ \rho(\per{1A0},A) = \lim_{n\to\infty}\frac{n}{n(|A|+2)}
                     = \frac{1}{|A|+2} < \frac12
  = \lim_{n\to\infty} \frac{n|A|/2}{n|A|} = \rho(\map(\per{1A0}),A) .$$
In the second case we choose $x:=\per{0A1}$ to come to a similar contradiction.
\qed

\n{\bf Proof} of Theorem~\ref{t:conf-conv} in the case $v=M=1$
follows from Lemmata~\ref{l:conv-low}-\ref{l:per-str}. \qed

\n{\bf Proof} of Theorem~\ref{t:flux-conv} in the case $v=M=1$.
We have the following identity:
$\rho_\pm(x,1)=\rho_\pm(x,10)+\rho_\pm(x,11)$. If
$\rho_+(x,1)\le1/2$ then $\map^tx\to \Free$
and $\rho_+(\map^tx,11)\to0$, thus
$\Phi_\pm(\map^tx)=\rho_\pm(\map^tx,10)
 = \rho_\pm(\map^tx,1)-\rho_\pm(\map^tx,11)
   \to\rho_\pm(\map^tx,1)=\rho_\pm(x,1)$.
The situation $\rho_-(x,1)\ge1/2$ can be reduced to the previous one by
going to the dual configuration.

Consider now the case $\rho_-(x,1)<1/2<\rho_+(x,1)$. By definition
there exists a sequence of pairs of positive integers
$n'_{i},m'_{i}\to\infty$ such that
$\rho(x[-n'_i,m'_i],1)\toas{i\to\infty}\rho_-(x,1)<1/2$. For each
$i$ we choose integers $n_{i}\ge n'_{i}$, $m_{i}\ge m'_{i}$ to be
the smallest integers satisfying the condition that $-n_{i}-1$ is
the ending point and $m_{i}+1$ is the starting point of some
nonoverlapping minimal intervals of the configuration $x$. If there
are no more nonoverlapping minimal intervals in the considered
direction or the segment $x[-n'_i,m'_i]$ intersects with an infinitely
long minimal interval we set $n_i:=n_i'$ or $m_i:=m_i'$ respectively,
depending on the direction where this event occurs. Clearly, we have
$\rho(x[-n_i,m_i],1)\le\rho(x[-n'_i,m'_i],1)$ and thus
$\rho(x[-n_i,m_i],1)\toas{i\to\infty}\rho_-(x,1)<1/2$. By the
definition of minimal intervals after $t_i:=(n_i+m_i)/2+1$
iterations all clusters of particles inside of the segment
$x[-n_i,m_i]$ will disappear and all particles will become free.
Therefore we can again apply the same argument as in the case
$\rho_+(x,1)\le1/2$ and obtain the relation for the lower limit of
the flux. To obtain the relation for the upper limit one should
consider the dual configuration. \qed

\section{Dynamics of fast particles $(T_v,X)$} \label{s:fast}

Note that the analysis of dynamics of the slow particles model
$(\map,X)$ is divided logically into two parts: first, we study
low density initial configurations $x\in X$ with $\rho_+(x,1)\le1/2$,
and then for high density configurations $x\in X$ with
$\rho_-(x,1)>1/2$ we pass to the dual ones using the property
that $\rho_+(\c{x},1)\le1-\rho_-(x,1)<1/2$ and argue that the
dual map $\c\map\equiv\map_{-1}$ has exactly the same asymptotic
properties as $\map$. The problem with the fast particles model
$(\map_v,X)$ is that the dual map map $\c{\map_v}\ne\map_{-v}$
in this case, and, in fact, has a very nontrivial dynamics.
Namely, $\c{\map_{-v}}$ corresponds to the situation, known
in physical literature (in the case $v=2$) as a traffic model
with `smart drivers', who anticipating the motion of at most
$v$ cars ahead, may move to an occupied site ahead of it
with the maximal velocity $1$. Example for the case $v=2$:
$\per{01110}\toas{\map_2}\per{01011}$.

Therefore since we are unable to study directly the dual map in
this case and according to the entire ideology of this paper, we
elaborated a reduction to the main case $v=M=1$ based on the
following consideration. Note that under the action of the map
$\map$ on $x\in X$ each pair $10$ goes to $01$ (i.e. the position
of a particle and a hole are exchanged). Therefore $\map$
is equivalent to the substitution rule $10\to01$. To apply this
idea to the case of $\map_v$ we introduce an alphabet
$\cA_v:=\{0_1,0_2,\dots,0_v,1\}$ with $v+1$ symbols and a map
$C_v:X\to\cA_v^{\IZ}\equiv X_v$ defined as follows: for each
segment $x[i,i+n+1]=1\underbrace{0\dots0}_{n}1$, we set
$C_vx[i,i+n+1]:=1\underbrace{0_v\dots0_v}_{\intp{n/v}}0_{n-\intp{n/v}v}1$.
If $n-\intp{n/v}v=0$ we shall drop the last element in
$C_vx[i,i+n+1]$. It remains to define the action of $C_v$ on
`tails' of $x$ consisting of only zeros, which we set according
to the following rules:
$\dots0001\dots\toas{C_v}\dots0_v0_v0_v1\dots$ and
$\dots1000\dots\toas{C_v}\dots10_v0_v0_v\dots$.

Now we are ready to define the substitution map $S_v:X_v\to X_v$
acting in the set $X_v$ according to the set of $v$ substitution
rules $10_i\to0_i1$ for $0<i\le v$, which generalizes the
substitution rule for the slow particles dynamics for the case of
$v$ different types of holes.

To study the life-time of clusters of particles in configurations
$x\in X_v$ we introduce also a new map $\tilde\map:=C_vC_v^{-1}S_v C_v$.

\begin{lemma}\label{l:fast-equiv} $\map_v=C_v^{-1} S_v C_v$, and
$\map_v^n=C_v^{-1}\tilde\map^n$ for any $n\in\IZ_+$.
\end{lemma}

\proof Straightforward. Note only that the map $C_vC_v^{-1}$ needs
not to be identical, example:
$$ C_vC_v^{-1}\per{10_{v-1}0_11}=\per{10_v1} .$$ \qed

Observe that the map $\tilde\map$ acts on $X_v$ in exactly the same way
as $\map$ acts on the space of binary sequences, namely $\tilde\map$
moves each particle by one position forward if there is no particle
there or the particle preserves its position otherwise. So the only
difference is that now we have $v$ different types of zeros, instead
of the only one type in the case $v=1$.

Therefore to study the life-time of a cluster of particles we
apply a similar machinery as in the case of the slow particles.
Denote
$$ {\rm Ind}(a):=\function{-v &\mbox{if } a=1 \\
                            i &\mbox{if } a=0_i}, \quad
   I(A,i) := \max\{k<i: ~ \sum_{j=k+1}^i {\rm Ind}(A_j)<0\}, \CR
   \Omega^{n} := \{A\in\cA_M^n: ~ A_n=1, ~ I(A,n)=1\} ,$$
where $A_j$ is (as usual) the $j$-th element of the word $A$.
Consider a map $\Gamma$ defined on words of length $n\in\IZ_+$
as follows: $(\Gamma A)_i:=(\tilde\map A)_{i+1}$ for all $i=1,2,\dots,n-1$.

\begin{lemma}\label{l:time-leve-v}
$\Gamma:\Omega^{n}\to\Omega^{n-2-\xi}$, where $0\le\xi<n-1$. The
life-time of the cluster of particles in the end of a word
$A\in\Omega^n$ is equal to $(\rho(A,1)\cdot|A|-1)$.
\end{lemma}

\n{\bf Proof} follows from the same argument as the one of
Lemma~\ref{l:life-time}. The only difference is that due to the
action of $\Gamma$ the number of elements in $\Gamma A$ may
become smaller than $|A|-2$, since the action of $C_vC_v^{-1}$
may decrease the number of $0_i$. On the other hand, during the one
iteration of the map $\Gamma$ only one elemement $1$ disappears
from (the right hand side) of $A$, i.e.
$\rho(A,1)\cdot|A| = \rho(\Gamma A,1)\cdot|\Gamma A| + 1$. Therefore
the number of iterations needed for the cluster of particles in the
end of the word $A$ to disappear is equal to the number of ones in
the word $A$ minus one. \qed

\begin{center}
\begin{tabular}{rccccccccc}
$A=$        & $0_2$& $0_2$& $0_1$&   $1$& $0_1$&   $1$& $0_1$& $1$& $1$ \\
$\Gamma^1A=$&      &      & $0_2$& $0_2$&   $1$& $0_1$&   $1$& $1$&     \\
$\Gamma^2A=$&      &      & $0_1$& $0_2$&   $1$&   $1$&      &     \\
$\Gamma^3A=$&      &      &      & $0_2$&   $1$&      &      &     \\
\end{tabular}\par
\smallskip
Example of the action of $\Gamma^t$ on $\Omega^{9}$ with $v=2$.
\end{center}

\bigskip

\n{\bf Proof} of
Theorems~\ref{t:density-pres},\ref{t:conf-conv},\ref{t:flux-conv}
for the case $v>1, M=1$ follows immediately from
Lemmata~\ref{l:fast-equiv},\ref{l:time-leve-v} and the reduction to
the case $v=1$ obtained there. \qed

Consider now a special case of {\em superfast} particles corresponding
to the choice of maximal velocity $v=\infty$. Denote
$$ X^{(\infty)}:=\{x\in X: ~~ \forall n\in\IZ ~~ \exists m,m'>|n|:
                         ~ x_m=1, ~ x_{-m'}=0\} ,$$
i.e. the set of binary configurations having no infinitely long right
`tails' of zeros or left `tails' of ones. Then the maps
$\map_{\infty}, \c{\map_{\infty}}: X^{(\infty)}\to X^{(\infty)}$
are well defined. The substitution rule
$1\underbrace{0\dots0}_{i}1\to10_{i}1$ $\forall i\in\IZ_+$ maps
$X^{(\infty)}\to X_\infty$. Strictly speaking, the latter has an
infinite alphabet, however all arguments applied in the case of
finite $v$ work as well. Moreover here the situation is even simpler,
because between each pair of consecutive ones there is only one
zero with a certain finite index: $\dots 1 0_i 1 0_j \dots$. Thus
the dynamics of $(\map_{\infty},X_\infty)$ is equivalent to the
dynamics of free particles, which gives the flux
$\Phi_\pm(x)=1-\rho_\pm(x,1)$.

\section{Dynamics of multi lane flows $(\map_{v,M},X_M)$.\\
         Reduction of $\map_{1,M}$ to the direct product of $M$ maps $\map$}
\label{s:dir-prod}

The model of a muti lane flow of slow particles on a finite
lattice has been introduced in \cite{NT} and generalized for the
case of an infinite lattice $\IZ$ in \cite{Bl-jam}, where statistical
properties of regular initial configurations have been obtained.
However the approach used in \cite{Bl-jam} does not allow to study
the dynamics of general initial configurations, which we shall
consider in this Section using a completely different method.

Our first aim is to redistribute a configuration $x\in X_M$ into $M$
binary configurations $\{x^{(j)}\in X=X_1\}_{j=1}^M$, such that
$\map_{1,M}^t x=\sum_j \map^tx^{(j)}$ for all $t\in\IZ_+\cup\{0\}$,
where the notation $x=\sum_j x^{(j)}$ means that $x_i=\sum_j x^{(j)}_i$
for each $i\in\IZ$. To solve this problem we introduce a 
{\em sawtooth redirection} $S_l:X_M\to(X_1)^M$ with 
$S_l x=\{x^{(j)}\}_{j=1}^M$ of a configuration $x\in X_M$ to a
collection of binary configurations $\{x^{(j)}\}_{j=1}^M$ with 
the starting point at site $l\in\IZ$: 
$$ x_i^{(j)} := \function{
   1 &\mbox{if } i\ge l \mbox{ and }
       j\in (\bigoplus_{k=l}^{i-1}x_k, \bigoplus_{k=l}^{i}x_k] \\
   1 &\mbox{if } i< l \mbox{ and }
       j\in (\bigoplus_{k=i}^{l-1}(-x_k), \bigoplus_{k=i+1}^{l-1}(-x_k)] \\
   0 &\mbox{otherwise,}} $$
where $a\oplus b:=(a+b-1)({\rm mod}{M})+1$ and
$\bigoplus_{i=n}^{m}x_{i}:=x_{n}\oplus \dots\oplus x_{m}$.
In other words, for the configuration $x\in X_{M}$ we construct a
beinfinite `staircase' starting from the site $l$ with the $i$-th
stair of height $x_{i}$ and then redistribute the result modulo
$M$ (preserving the site number) among $M$ binary configurations
$\{x^{(j)}\}_{j=1}^M$.

With some abuse of notation we shall refer to
$(S_{l}x)^{(j)}\equiv x^{(j)}$ as the $j$-th lane of
$S_{l}x=\{x^{(j)}\}_{j=1}^M$ and denote the
action of the direct product of maps $\map_v$ applied at $S_{l}x$ as
$\map_v S_{l}x:=\{\map_v x^{(j)}\}_{j=1}^M$.

Example for the case $v=1,M=3$ and the starting site $l$ defined by
the relation $x[l-1,l+1]=211$:
$$             S_l(\dots11211221\dots) ~
 = \begin{array}{l}\dots10100110\dots\\
                   \dots00101010\dots\\
                   \dots01010101\dots\end{array} \toas{\map_{1,3}}
   \begin{array}{l}\dots01010101\dots\\
               \dots\ast0010101\dots\\
               \dots\ast010101\ast\dots\end{array}
 = S_{l+1}(\dots\ast112121\ast\dots),$$
where the unknown positions are marked by $\ast$.

Symbolically the sawtooth redirection is shown in
Fig.~\ref{fig-sawtooth}(b) by curvilinear lines corresponding to
sawtooth rows of ones, open circles mark the intersections of
these lines with the `lanes' $j,j'$, i.e. the positions where
$x^{(j)}$ or $x^{(j')}$ are equal to $1$ (all other positions on
these lanes are occupied by zeros).
%%%%%%%%%%%%%%%%%%%%%%%%%%%%%%%%%%%%%%%%%%%%%%%%%%%%%%%%%%%%%%%%%%
%% Sawtooth redirection
\Bfig(300,100) %(150,100)
      {\footnotesize{
      \put(150,0){\put(70,-10){(b)}  %Fig.(b)
       \bline(0,0)(1,0)(150)   \bline(0,0)(0,1)(100)
       \bline(0,100)(1,0)(150) \bline(150,0)(0,1)(100)
       \bezier{60}( 7,20)(100,20)(145,20) \put(2,18){$j$}
       \bezier{60}( 7,80)(100,80)(145,80) \put(2,78){$j'$}
       \bezier{200}(10,10)(20,60)(35,90)
              \put(12,20){\circle3} \put(31,80){\circle3}
       \bezier{200}(35,10)(70,20)(70,90)
              \put(52,20){\circle3} \put(70,80){\circle3}
       \put(85,50){$\dots$}
       \bezier{200}(100,10)(120,80)(140,90)
              \put(103,20){\circle3} \put(129,80){\circle3}
       \bezier{35}(10,10)(10,50)(10,90)
       \bezier{35}(35,10)(35,50)(35,90) \put(33,3){$i_{-}$}
       \bezier{35}(70,10)(70,50)(70,90) \put(68,3){$i_{+}$}
       \bezier{35}(100,10)(100,50)(100,90) \bezier{35}(140,10)(140,50)(140,90)
      }
      \put(0,0){\put(70,-10){(a)}  %Fig.(a)
       \def\BOX#1{\bline( 0,0)(1,0)(10) \bline(0, 0)(0,1)(10)
           \bline(10,0)(0,1)(10) \bline(0,10)(1,0)(10) \put(2,3){#1}}
       \bline(0,0)(1,0)(150)   \bline(0,0)(0,1)(100)
       \bline(0,100)(1,0)(150) \bline(150,0)(0,1)(100)
       %\put(0,0){\BOXX{-3}{10}}\put(10,0){\BOXX{2}{10}}
       \put(35,40){\BOX0}                    \put(37,33){l}
       \put(45,50){\BOX1} \put(45,60){\BOX2} \put(47,33){l+1}
       \put(55,70){\BOX3}                    %\put(57,43){l+2}
       \put(25,30){\BOX{-1}} \put(25,20){\BOX{-2}} \put(25,43){l-1}
       \bline(15,20)(1,0)(10) \put( 5,10){\BOX{-3}}
       \put(80,40){\BOX{-3}} \bline(90,50)(1,0)(10)
       \put(100,50){\BOX{-2}} \put(100,60){\BOX{-1}}
       \put(110,40){\BOX0}                           \put(112,33){l}
       \put(120,50){\BOX1}\put(120,60){\BOX2}
       \put(130,40){\BOX3}
       \put(60,45){\vector(1,0){15}} \put(61,47){$S_l$}}
      }}
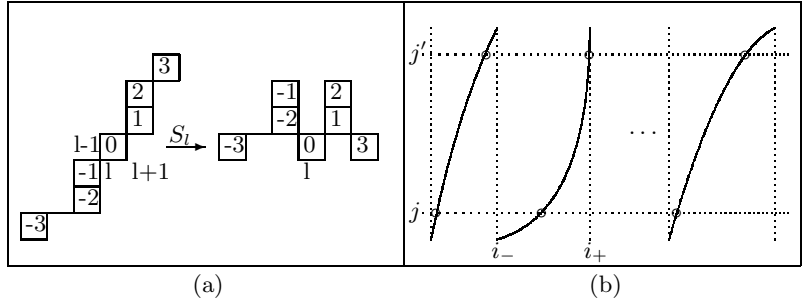
{`Sawtooth redirection':
(a) action of the map $S_l$ on individual particles marked by squares with
their relative numbers ($M=3$);
(b) symbolical representation of $S_{l}x$, $[i_-,i_+]$ -- one of intervals
of monotonicity.%
\label{fig-sawtooth}}
%%%%%%%%%%%%%%%%%%%%%%%%%%%%%%%%%%%%%%%%%%%%%%%%%%%%%%%%%%%%%%%%%

\begin{theorem}\label{t:sawtooth} For any $x\in X_M$ and $l\in\IZ$
and $S_l(x)\equiv\{x^{(j)}\}_{j=1}^M$ we have\par%
(a) $x=\sum_j (S_{l}x)^{(j)}$,\par%
(b) $|\rho(x^{(j)}[n+1,n+k],1) - \rho(x^{(j')}[n+1,n+k],1)|\le1/k$
    ~ $\forall j,j'\in\{1,\dots,M\}, ~ n\in\IZ$ and $k\in\IZ_+$,\par%
(c) $S_{l+k}x=\{(S_{l}x)^{(j\oplus k)}\}_{j=1}^M$ ~ $\forall k\in\IZ_+$,\par%
(d) $\forall v\ge1$ we have $\map_v S_{l}x=S_{l+\xi}x'$
    for some $\xi\in\{0,1,\dots,v\}$ and $x'\in X_M$
    which doesn't depend on $l$,\par%
(e) $\map S_{l}x= S_{l+\xi}(\map_{1,M}x)$ for some $\xi\in\{0,1\}$.
\end{theorem}

\proof The statement (a) follows immediately from the definition of
the sawtooth redirection, because during the redirection each
particle preserves its position $i$.

The property (b) is equivalent to the assumption that
$$ \left|\sum_{i=1}^k x_{n+i}^{(j)}
       - \sum_{i=1}^k x_{n+i}^{(j')}\right| \le 1 ,$$
i.e. that the number of particles in the same segment of different
`lanes' $j,j'$ can differ at most by $1$. According to the
`sawtooth redirection' for any given finite segment of integers
$n+1,n+2,\dots,n+k$ the number of intersections of the curvilinear
lines in Fig.~\ref{fig-sawtooth}(b) with the horizontal line at
height $j$ differs from number of intersections with the horizontal
line at height $j'$ at most by one.
This immediately yields the property (b).

The collection of binary configurations $S_l(x)$ has a row of ones
at site $i$ of height $k$ if and only if $x_i=k$, and the change
of the starting point $l$ of the redirection only changes
cyclically the starting point $1$ of the enumeration of lanes
$x^{(j)}$. This proves the property (c).

Observe now that the definition of $S_l(x)$ is equivalent to the
existence of a partition of $\IZ$ into segments $[i_-,i_+]$ such that
$x_{i_-}^{(1)}=1$, $x_{i_+}^{(M)}=1$ (except for the most left segment 
where $i_-=-\infty$ and the most right one where $i_+=\infty$) and 
for any $1<j<M$ there exists the only one $i\in[i_-,i_+]$ such that
$x_{i}^{(j)}=1$. Indeed, according to the `sawtooth redirection' the
curvilinear lines in Fig.~\ref{fig-sawtooth}(b) have the property that the
intersection with the horizontal line at height $j$ occurs not earlier
than with the horizontal line at height $j'>j$ (the curvilinear lines
may have vertical segments). To simplify the notation we shall say that
$S_l(x)$ is monotonous on $[i_-,i_+]$.

Consider the interval of monotonicuty $[i_{-},i_{+}]$ which starts
from $l$, i.e. $i_{-}=l$. We set $\xi$ to be equal to the minimum
of the number of not occupied positions in $x^{(1)}$ ahead of the
site $i_-$ (which is occupied by 1). Then under the action of
$\map_v$ the particle at the site $i_-$ of the 1-st lane moves by
$\xi$ positions to the right. Observe that all particles on the
other lanes in the segment $[i_-,i_+]$ have at least $\xi$ not
occupied positions ahead of them, and therefore all these
particles will move at least $\xi$ positions to the right. Thus to
prove that the monotonicity is preserved it is enough to note that
the particle on the lane $M$ cannot move further to the right
than the first particle on the first lane of the next interval of
monotonicity. Indeed, the latter is a trivial consequence of the
definition of intervals of monotonicity. This finishes the proof
of the statement (d) except the last part, which follows from the
statement (c). %\met

To prove the statement (e), observe that by the definition of the
map $\map_{1,M}$ (see Section~\ref{s:intro}) a particle at the
site $i$ of the lane $j$ can switch to the lane $j'$ if and only
if $x^{(j)}[i,i+1]=11$ and $x^{(j')}[i,i+1]=00$, which
contradicts to the definition of the intervals of monotonicity.
Therefore under the sawtooth redirection no particle in $S_l(x)$
will change its lane. \qed

\begin{corollary} The sawtooth redirection gives a simple constructive
way to rearrange vehicles in a multi lane traffic flow between lanes
(preserving their positions in the flow) to achieve the maximal
available flux.
\end{corollary}

According to Theorem~\ref{t:sawtooth},(d) the map
$x\to\sum_{j}\map_{v}(S_{l}x)^{(j)}$ is well defined as a map
from $X_M$ into itself and does not depend on the choice of the
starting site $l\in\IZ$. Moreover, it can be shown that this
formula coincides with (\ref{def:fast-one}) in the case $v=1$,
and it clearly coincides with $\map_v$ in the case $M=1$.
Therefore we use this relation as a definition of the dynamics
of a general multi lane flow in the case $v,M>1$, namely we set
$\map_{v,M}x:=\sum_{j}\map_{v}(S_{0}x)^{(j)}$.

\bigskip

\n{\bf Proof} of Theorem~\ref{t:density-pres} for the case
$v,M\ge1$ and $A\subset X$ follows now from the sawtooth
redirection, which gives the reduction to the one-lane case. It
remains to show that the statistics of more general words
$A\subset X_M$ with $|A|=1$ might be not preserved under
dynamics. The reason of this is that if $M>1$ the multiplicities
might be not preserved. Indeed, let
$a\in\cA_{M}\setminus\{0,1\}$. Then
$\rho(\per{a(M-a+1)0},a)=\frac13(1+\intp{(M-a+1)/a}+0)$, while%
\bea{ \rho(\map_{1,M}\per{a(M-a+1)0},a)
  \a= \rho(\per{1(a-1)(M-a+1)},a) \\
  \a= \frac13(0+0+\intp{(M-a+1)/a}+0) < \rho(\per{a(M-a+1)0},a) .}%
\qed

\n{\bf Proof} of Theorems~\ref{t:conf-conv},\ref{t:flux-conv}
for the case $v,M>1$. Consider a configuration $x\in X_M$.
According to Theorem~\ref{t:sawtooth},(b) for
$S_0x\equiv\{x^{(j)}\}_{i=1}^M$ we have $\forall n,m\in\IZ_+$ that
$$ |\rho(x^{(j)}[-n,m],1) - \rho(x^{(j')}[-n,m],1)|\le\frac1{m+n+1} .$$
Thus going to the limit as $n,m\to\infty$ and using
Theorem~\ref{t:sawtooth},(a) we get
$\rho_\pm(x^{(j)},1)=\frac1M \rho_\pm(x,1)$ for each $j\in\{1,\dots,M\}$.
Therefore the application of the results obtained in
Sections~\ref{s:slow},\ref{s:fast} in the case of one-lane flows (i.e.
in the case of the map $\map_v$) proves the statements under question. \qed

\section{Rate of convergence: (space) periodic, regular, and typical \\
         initial configurations} \label{s:rate}

In this section we study the rate of convergence of various
statistics to the corresponding limit values whose existence have
been established in Theorems~\ref{t:flux-conv},\ref{t:conf-conv}.
Since we have shown that the analysis of $\map_{v,M}$ in all
cases can be reduced to the case of $v=M=1$, we consider in this
section only the dynamics of slow particles and shall consider the
proof of Theorem~\ref{t:typ-conv} only for this case.

We start with periodic in space configurations. Clearly each
$n$-periodic in space configuration $x\in X$ can be represented
in the form $x:=\per{A}$ with a binary word $A$ of length $|A|=n$.

\begin{lemma} Let $x:=\per{A}$ with $|A|=n$, then
$\rho_\pm(x,a)=\rho(A,a)$ for $a\in\{0,1\}$. The space of
$n$-periodic in space configurations is invariant under
the action of the map $\map$ and after at most $\intp{n/2}+1$
iterations any configuration from this space belongs to
$\Free\cup\c\Free$. \end{lemma}

\proof Straightforward. \qed

The only nontrivial question related to this space is the length
of the transient period for a given configuration $x\in X$.

\begin{lemma} Let $x:=\per{A}$ with $|A|=n$ and $\rho(A,1)\le1/2$
and let $B$ be the longest minimal word in $A$. Then the length
of the transient perion is equal to $|B|/2-1$.
\end{lemma}

\proof This is an immediate consequence of Lemma~\ref{l:life-time}.
\qed

Consider now a generalization of the space of periodic in space
configurations -- the space of regular configurations, proposed
in \cite{Bl-var,Bl-jam}. This space is defined as follows:
$$ \Reg(r,\psi) := \{x\in X: ~ |\rho(x[-n,m],1) - r|\le\psi(n+m)
                     ~~ \forall n,m\in\IZ_+\} ,$$
where  $r\in[0,1]$ is a constant, and the nonnegative function
$\psi(n)\toas{n\to\infty}0$ is assumed to be a strictly decreasing.

\begin{lemma} Let $x\in\Reg(r,\psi)$, then the density $\rho(x,1):=r$
is well defined, and if $\ne1/2$ then there is a constant $\tau$ such
that the life-time of any cluster of particles in $x$ does not exceed $\tau$.
\end{lemma}

\proof First, observe that if $x\in\Reg(r,\psi)$ then we have
$$ \rho_\pm(x,1)=\blim{\pm}{n,m\to\infty}\rho(x[-n,m],1)=r ,$$ 
and thus $\rho(x,1)=r$ is well defined. Assume now that $x\in\Reg(r,\psi)$
with $r<1/2$, then, since $\psi(n)\toas{n\to\infty}0$, we deduce that
there is a positive integer $\tau$ such that $\phi(\tau) < 1/2 - r$.
Then we have $\rho(x[n,n+\tau],1)<1/2$ for any $n\in\IZ$, which yields
the claim of the lemma due to Lemma~\ref{l:conv-low}. The case $r>1/2$
follows from the same argument but applying Lemma~\ref{l:conv-high}
instead. \qed

Now we proceed to study more general initial configurations.

\begin{lemma}\label{l:conv-simple}
Let $x\in X$ satisfies the assumption that there exists a number
$\gamma\in(0,1)$ such that $\forall n\in\IZ_{+}$ and any word
$A\subseteq x[-n,n]$ with $|A|>2\gamma n$ we have
$\rho(A,1)\le1/2$. Then
$\dist(\map^{t}x,\Free)\le2^{-t/\gamma +1}$ for any $t\in\IZ_+$.
If $\c{x}$ satisfies the same assumption, then we have
$\dist(\map^{t}x,\c\Free)\le2^{-t/\gamma +1}$.
\end{lemma}

\proof Consider only those $n\in\IZ_{+}$ for which the largest minimal
words containing a cluster of particles in $x[-n,n]$ also belong to
$x[-n,n]$. By the assumption of Lemma the length of the largest minimal
interval containing in the segment $x[-n,n]$ does not exceed $2\gamma n$.
Therefore the corresponding clusters of particles with disappear after
at most $\gamma n$ iterations, and thus for all sufficiently
large $t\in\IZ_{+}$ all particles in the segment
$\map^{t}x[-t/\gamma, t/\gamma]$ will become free. Thus the closest to
the origin nonfree particle can appear not earlier as at site $t/\gamma$,
which gives the desired estimate of the rate of convergence. The second
statement follows from the same argument applied to the dual map.
\qed

\begin{lemma}\label{l:conv-flux}
Let $x\in X$ satisfies the same assumption as in
Lemma~\ref{l:conv-simple}, then
$$ \blim+{n\to\infty}\frac1{2n}\sum_{i=-n}^nV(\map^nx,1) 
 = F_{1,1}(\rho(x,1)) .$$
\end{lemma}

\proof Observe that $\frac1{2n}\sum_{i=-n}^nV(x,1)=\rho(x[-n,n],10)$.
Applying the same argument as in the proof of Lemma~\ref{l:conv-simple}
we see that after $n$ iterations the segment $\map^{n}x[-n/\gamma, n/\gamma]$
contains only free particles. Therefore $\rho(x[-n,n],10)=\rho(x[-n,n],1)$,
which yields the desired equality. \qed

\begin{corollary} The statements of
Lemmata~\ref{l:conv-simple},\ref{l:conv-flux} remain valid if instead of
$\forall n\in\IZ_{+}$ we assume that $n$ belongs to the subset of
$\IZ_{+}$ of density 1.
\end{corollary}

\begin{lemma}\label{l:typ-conf} $\forall \gamma\in(0,1)$ for
$\mu_p$-a.a. configurations $x\in X$ the set of $n\in\IZ_{+}$,
for which any word $A\subseteq x[-n,n]$ with $|A|>2\gamma n$
satisfies the inequality $\rho(A,1)\le1/2$, has the density 1.
\end{lemma}

\proof\hskip-9pt\footnote{The idea of this construction, based on 
the large deviation principle, was proposed by A.~Puhal'skii.}
Let $\{x_i\}_{-\infty}^{\infty}$ be a Bernoulli sequence with
the density $\cP(x_i=1)=p<1/2$ for all $i\in\IZ$. Introduce
a sequence of functions
$y_n(\tau):=\frac1{2n+1}\sum_{i=-n}^{-n+\intp{2n\tau}}x_i$
depending on a real variable $\tau\in[\gamma,1]$, and consider 
a functional
$$ \phi(y(\tau)) := \sup_{\tau\in[0,1-\gamma]}
                    \sup_{\gamma\le s\le1-\tau}\frac1s(y(\tau+s)-y(\tau)) $$
defined in Skorohod space of functions $y(\tau)$.
Then the quantity under question is the probability
$\cP(\phi(y_n(\tau))\le1/2 ~~ \forall\tau\in[\gamma,1])$.
Since $y_n(\tau)$ converges in probability for a given $\tau$ to
$\tilde y(\tau):=p\tau$ and the functional $\phi$ is continuous,
$\phi(y_n(\tau))$ converges to $\phi(\tilde y(\tau))$ (functional law
of large numbers). Thus we have
$$ \cP(\phi(y_n(\tau))\le1/2 ~~ \forall\tau\in[\gamma,1])
   \to \cP(\phi(\tilde y(1))\le 1/2) = 1 ,$$
where the rate of convergence
$(\cP(\phi(y_n(1))>1/2))^{1/n} \toas{n\to\infty} \sqrt{2p(1-p)}$
follows by the combination of the large deviation principle for the
functions $y_n(\tau)$ and the contraction principle (see, e.g., \cite{DZ}).
\qed

\begin{corollary} Results of
Lemmata~\ref{l:conv-simple},\ref{l:conv-flux},\ref{l:typ-conf}
prove Theorem~\ref{t:typ-conv} in the case $v=M=1$.
\end{corollary}

\section{Dynamics of measures and chaoticity}
\label{s:measure}

In this section we shall study the action of the map $\map_{v,M}$
in the space $\cM(X_{M})$ of probabilistic measures on $X_{M}$. 
This action is defined as follows:
$\map_{v,M}\mu(Y):=\mu(\map_{v,M}^{-1}Y)$ for a measure
$\mu\in\cM(X_{M})$ and a measurable subset $Y\subseteq X_{M}$. A
measure $\mu\in\cM(X_{M})$ is called translation invariant if
it is invariant with respect to the action of the shift map
$\sigma:X_{M}\to X_{M}$.

\begin{lemma} If $\mu\in\cM(X_{M})$ is translation invariant
then this property holds for $\map_{v,M}^{t}\mu ~ \forall
t\in\IZ_{+}$.
\end{lemma}

\proof We have $\map_{v,M}^{t}\mu(Y) = \mu(\map_{v,M}^{-t}Y) =
 \mu(\sigma\map_{v,M}^{-t}Y) =
 \mu(\map_{v,M}^{-t}\sigma Y) = \map_{v,M}^{t}\mu(\sigma Y)$.
\qed

One might expect that under the action of the map $\map_{v,M}$
any translation invariant measure should converge to a Bernoulli 
one. Indeed,%
\bea{\map\mu_p(x\in X: ~ x_0=1)
 \a= \mu_p(x\in X: ~ x[0,1]=11) + \mu_p(x\in X: ~ x[-1,0]=10) \\
 \a= \mu_p(x\in X: ~ x[0,1]=11) + \mu_p(x\in X: ~ x[0,1]=10)
   = \mu_p(x\in X: ~ x_0=1) .}%
On the other hand, the product structure is not preserved 
even in the case of the model of slow particles.

\begin{lemma} The measure $\map\mu_p$ is not a product one 
for any $0<p<1$.
\end{lemma}

\proof We have% 
\bea{\map\mu_p(x\in X: ~ x[0,1]=11) 
 \a= \mu_p(x\in X: ~ x[0,2]=111) + \mu_p(x\in X: ~ x[-1,2]=1011) \\
 \a= p^3+p^3(1-p)=p^3(2-p) \ne p^2 = \mu_p(x\in X: ~ x[0,1]=11) .}%
Thus the measure $\map\mu_p$ does not have the product structure. \qed

It is of interest that in the case $v>1$ even the average value 
$\mu_p(x\in X: ~ x_0=1)$ is not preserved under dynamics. Indeed,
\bea{\map_v\mu_p(x\in X: ~ x_0=1)
 \a= \mu_p(x\in X: ~ x[0,1]=11) + \sum_{i=1}^v\mu_p(x\in X: ~ x[-1,0]=10_i) \\
 \a= p^2 + p(1-p) + \dots + p(1-p)^v = p + p\sum_{i=2}^v(1-p)^i \\
 \a> p = \mu_p(x\in X: ~ x_0=1) .}%

Note that in the case of the slow particles model $(\map_{1,1},X)$ 
some results about the set of $\map_{1,1}$-invariant measures and 
mathematical expectations of the limit flux with respect to them 
were studied in \cite{BKNS}. 
\smallskip

In \cite{Bl-var} it has been proven that the dynamical system
$(\map_{1,M},X_M)$ is chaotic in the sense that its topological
entropy is positive. Moreover this paper gives an asymptotically
exact (as $M\to\infty$) representation for the entropy. The extension
of this result to the case $(\map_{v,M},X_M)$ with $v>1$ is
straightforward.

\section{Passive tracer in the 1-lane flow of fast particles}\label{s:tracer}

Let $\map_v^tx$, $v\ge1$ describes the 1-lane flow of particles and
let at time $t$ the passive tracer occupies the position $i$. Then before
the next time step of the model of the flow the tracer moves in its
chosen direction to the closest (in this direction) position of a
particle of the configuration $\map_v^tx$. For example, if the going
forward tracer occupies the position 2 and the closest particle
in this direction occupies the position 5, then the tracer moves
to the position 5. Then the next iteration of the flow occurs,
the tracer moves to its new position, etc.

To be precise let us fixed a configuration $x\in X$ with
$\rho_-(x,1)>0$ and introduce the maps $\tau_{x}^{\pm}:\IZ\to\IZ$
defined as follows:
$$ \tau_{x}^{+}i := \min\{j: \; i<j, \; x_j=1\}, \quad %\CR
   \tau_{x}^{-}i := \max\{j: \; i>j, \; x_j=1\} .$$
Then the simultaneous dynamics of the configuration of particles
(describing the flow) and the tracer is defined by the skew product of
two maps -- the map $\map_v$ and one of the maps $\tau_{\cdot}^{\pm}$, i.e.
$$ (x,i) \to {\cal T}_{\pm}(x,i) := (\map_v x, \tau_{x}^{\pm}i) ,$$
acting on the extended phase space $X\times\IZ$. The sign $+$ or
$-$ here corresponds to the motion along or against the flow.
We define the {\em average (in time) velocity} of the tracer $V(t,x)$
as $S(t)/t$, where $S(t)$ denotes the total distance covered by the
tracer (which starts at the site 0) up to the moment $t$ with the
positive sign if the tracer moves forward, and the negative sign otherwise.

\begin{theorem}\label{t:tracer}
Let $x\in\{x\in X: ~ \dist(\map_v^tx,\Free\cup\c\Free)\le2^{-t/\gamma+1}\}$ 
for all $t\in\IZ_+$ and some $0<\gamma<1$.
If $0<\rho_+(x,1)\le\frac1{v+1}$, then $V(t,x)\toas{t\to\infty}v$ if
the tracer moves along the flow (i.e. in the case ${\cal T}_+$), and
$\blim\pm{t\to\infty}V(t,x)=\frac{-1}{\rho_\pm(x,1)}+1$ in the opposite 
case. If $\rho_-(x,1)>1-\frac1{v+1}$ and the tracer moves
against the flow then $V(t,x)\toas{t\to\infty}-1$.
\end{theorem}

\n{\bf Remark}. The assumption about the initial configurations is satisfied
for $\mu_p$-a.a. $x\in X_M$ (see Theorem~\ref{t:typ-conv}).

\smallskip

\proof Since we assume that $\map_v^tx$ converges to the attractor
$\Free\cup\c\Free$ with the exponentially fast rate, then at the moment
$t\in\IZ_+$ we have an exponentially long (in $t$) interval of the
configuration $\map_v^tx$ consisting of only free particles or
free holes (depending on the density). As we shall show that $V(t,x)$
converges to a constant, then to study its value we can restrict
the analysis to the case $x\in\Free\cup\c\Free$.

Under the assumption $0<\rho_+(x,1)\le\frac1{v+1}$ we have
$\map_v^tx\toas{t\to\infty}\Free_v$.
In the case of ${\cal T}_+$ the tracer will run down one of the
particles and will follow it, but cannot outstrip. Indeed after
each iteration of the flow this free particle occurs exactly
$v$ positions aheed of the tracer. Thus $V(t,x)\toas{t\to\infty}v$.

Consider now the case when the tracer moves backward with respect 
to the flow. Then each time when the tracer encounters a particle, 
on the next time step this particle moves in the opposite direction 
and does not interfere with the movement of the tracer. We assume 
again that $x\in\Free_v$ and consider the case 
$0<\rho_+(x,1)\le\frac1{v+1}$. 
If on the spread of length $n$ there are $m$ particles,
i.e. $m$ obstacles for the tracer then the average velocity on
this segment is equal to $\frac{n-m}{m}$. Going to the limit as
$n\to\infty$ we obtain the desired estimate.

It remains to consider the case $\rho_-(x,1)>1-\frac1{v+1}$ and thus
$\map_v^tx\toas{t\to\infty}\c{\Free_v}$, i.e. to the flow where all
holes move at maximal velocity $-v$. Thus after each iteration the
tracer moves exactly by one position to the left (since it never can
encounter a hole), which gives the limit velocity $-1$. \qed

Observe that the motion against the flow is efficient only in the
case of low density of particles when $\rho_+(x,1)\le\frac1{v+1}$.
On the other hand, in the high density region in the case of the motion
along the flow and in the region $\frac1{v+1}<\rho_-(x,1)<1-\frac1{v+1}$
in the case of the motion against the flow the limit velocity of the
tracer depends not only on the densities, but also on the fine structure
of the configuration $x$. Moreover, this concerns also the case of
`untypical' initial configurations with $0<\rho_-(x,1)<1/2<\rho_+(x,1)$,
when there might be arbitrary long (even infinite) minimal words for
both particles and holes.

%\newpage
%\small

\end{document}